\documentclass[art]{amsart}

\usepackage{amssymb,latexsym}
\usepackage{amsfonts}
\usepackage{dsfont}
\usepackage[pdftex]{hyperref}
\usepackage{color}
\usepackage{graphicx}
\usepackage[all]{xy}
\usepackage{youngtab}

\theoremstyle{plain}
\newtheorem{thm}[equation]{Theorem}
\newtheorem{lem}[equation]{Lemma}

\newtheorem{prop}[equation]{Proposition}

\newtheorem{cor}[equation]{Corollary}
\newtheorem{rem}[equation]{Remark}
\newtheorem{defn}[equation]{Definition}
\numberwithin{equation}{section}

\newcommand{\C}{\mathbb C}
\newcommand{\fg}{\mathfrak{g}}
\newcommand{\fh}{\mathfrak{h}}
\newcommand{\fb}{\mathfrak{b}}

\newcommand{\fm}{\mathfrak{m}}

\newcommand{\fz}{\mathfrak{z}}

\renewcommand{\fz}{\mathfrak{z}}

\newcommand{\B}{\mathcal{B}}

\newcommand{\W}{\mathcal{W}}
\newcommand{\V}{\mathcal{V}}
\DeclareMathOperator{\Ad}{Ad}

\begin{document}

\title[The Betti numbers of regular Hessenberg varieties are palindromic]{The Betti numbers of regular Hessenberg varieties are palindromic}
\author{Martha Precup}
\address{Dept. of Mathematics, Northwestern University, Evanston, IL 60208}\email{mprecup@math.northwestern.edu}

%\keywords{Hessenberg varieties, palindromic}
%\subjclass[2010]{Primary 14L35, 14M15}

%\subjclass[2010]{Primary: 35Q53}
%
%

%
\maketitle
\begin{abstract}  
Recently Brosnan and Chow have proven a conjecture of Shareshian and Wachs describing a representation of the symmetric group on the cohomology of regular semisimple Hessenberg varieties for $GL_n(\C)$.  A key component of their argument is that the Betti numbers of regular Hessenberg varieties for $GL_n(\C)$ are palindromic.  In this paper, we extend this result to all reductive algebraic groups, proving that the Betti numbers of regular Hessenberg varieties are palindromic. 
\end{abstract}
 
 \section{Introduction}
 
The purpose of this paper is to show that the Betti numbers of regular Hessenberg varieties are palindromic.  The motivation to address this problem comes directly from the recent results in \cite{BC2} of Brosnan and Chow which analyze an action of the symmetric group on the equivariant cohomology of regular semisimple Hessenberg varieties for $GL_n(\C)$.  This action was defined by Tymoczko in \cite{T}.  Shareshian and Wachs conjectured that the character of this representation corresponds to the chromatic quasisymmetric function of an indifference graph (see \cite[\S 5]{SW}).  This conjecture is proven by Brosnan and Chow in \cite{BC2} and again using different methods by Guay-Paquet in \cite{G}.  

One key component of the monodromy argument used by Brosnan and Chow is the fact that regular Hessenberg varieties for $GL_n(\C)$ have palindromic Betti numbers (see \cite[Corollary 33]{BC}).  In this paper we extend this result to regular Hessenberg varieties for all reductive algebraic groups.  We hope these results will facilitate understanding of regular semisimple Hessenberg varieties for arbitrary reductive groups and the corresponding representation of the Weyl group on their cohomology. 

The main tools used here are the formulas given for the Betti numbers of Hessenberg varieties in \cite{P} and an analysis of subsets of roots of Weyl type by Sommers and Tymoczko in \cite{ST}.  The argument given is purely combinatorial, and recovers the result mentioned above for $GL_n(\C)$ due to Brosnan and Chow.  Their argument invokes results of Shareshian and Wachs proving that the chromatic quasisymmetric function of an indifference graph for a natural unit interval is palindromic (\cite[Corollary 4.6]{SW}).  Our arguments are Lie theoretic and do not use quasisymmetric functions, which are specific to the $GL_n(\C)$ case.  In addition, since regular Hessenberg varieties are singular in many cases, these results do not easily follow from Poincar{\'e} duality.  

%We note that the combinatorial methods used here are very different from those in \cite{SW} since the methods used in that paper involve quasisymmetric functions which are specific to the $GL_n(\C)$ case.  

The second section of this paper covers the background information and definitions needed for our arguments.   In the third, we use results from \cite{ST} and \cite{P} to analyze subsets of roots of Weyl type and prove our main result, which is Theorem \ref{main result} below.  We first prove a special case of the result; that the Betti numbers of regular nilpotent Hessenberg varieties are palindromic, stated in Proposition \ref{regular nilpotent} below.  This statement is part of an unpublished result of Peterson, appearing in \cite[Theorem 3]{BC}.  We provide a proof of Proposition \ref{regular nilpotent} which follows directly from the combinatorial results of Sommers and Tymoczko in \cite{ST}.

{\bf Acknowledgements.}  I would like to thank Patrick Brosnan and Sam Evens for helpful conversations and comments on this manuscript.  Part of this research was conducted during a Focused Research Group on Hodge Theory, Moduli and Representation Theory held at Texas A\&M University in January 2016.

%%%%%%%%%%%%%%%%%%%%%%%%%%%%%%%%%%

\section{Preliminaries}

Let $G$ be a reductive algebraic group with corresponding Lie algebra $\fg$ and Weyl group $W$.  Let $B$ be a Borel subgroup with Lie algebra $\fb$, $\fh$ be the Cartan subalgebra of $\fb$, and $U\subset B$ be the maximal unipotent subgroup.  Denote by $\Phi$ the associated root system, and by $\Phi^+$, $\Phi^-$ and $\Delta$ the subsets of positive, negative, and simple roots in $\Phi$, respectively.  Let $x \in \fg$ be a regular element.  Then $x$ can always be conjugated to an element of the form
\begin{eqnarray}\label{regelement}
x_J =n_J+s_J
\end{eqnarray}
 where $n_J=\sum_{\alpha\in J} E_{\alpha}$, $J\subseteq \Delta$ is a subset of simple roots, and $n_J$ is a regular nilpotent element in the Levi subalgebra $\fm_J:=\fz_{\fg}(s_J)$ of $\fg$.  Let $\Phi_J$ denote the root system generated by the simple roots $J\subseteq \Delta$ with positive roots $\Phi_J^+$ and negative roots $\Phi_J^-$.  Denote by $W_J$ the subgroup of $W$ generated by the simple reflections $s_{\alpha}$, $\alpha\in J$.  Every regular element corresponds to a subset $J$ of simple roots in this way so we will write $x_J$ for this regular element and assume it has been conjugated into the standard form given above.
  
\begin{defn}  A subspace $H\subseteq \fg$ is a {\it Hessenberg space} with respect to $\fb$ if $\fb\subset H$ and $H$ is a $\fb$-submodule of $\fg$, i.e., $[\fb,H]\subseteq H$.
\end{defn}
  
 Fix a Hessenberg space $H\subset \fg$ and let
\[
H= \fb \oplus \bigoplus_{\gamma\in \Phi_H^-} \fg_{\gamma}
\]
where $\Phi_H^-$ denotes the negative roots which correspond to negative root spaces in the Hessenberg space.  Set $\Phi_H=\Phi_H^- \cup \Phi^+$ and let $m_H=|\Phi_H^-|$.  
 
We denote by $\B=G/B$ the flag variety of $G$.  A Hessenberg variety is a subvariety of the flag variety defined as follows.

\begin{defn}  Fix $x\in \fg$ and a Hessenberg space $H$ with respect to $\fb$.  The Hessenberg variety associated to $x$ and $H$ is
		\[
				\B(x,H)=\{ gB \in \B : g^{-1}\cdot x\in H    \}
		\]
where $g\cdot x$ denotes the Adjoint action $\Ad(g)(x)$.  If $x$ is a regular element we say $\B(x,H)$ is a regular Hessenberg variety, if $n$ is a regular nilpotent element we say that $\B(n,H)$ is a regular nilpotent Hessenberg variety, and if $s$ is a regular semisimple element we say that $\B(s,H)$ is a regular semisimple Hessenberg variety.
\end{defn}

We denote the Schubert cell $BwB/B\subseteq \B$ by $C_w$, and let 
\[
N(w)=\{ \gamma\in \Phi^+: w(\gamma)\in \Phi^- \} \mbox{ and } N^-(w)=-N(w)=\{ \gamma\in \Phi^-: w(\gamma)\in \Phi^+ \}.		
\]
It is a well known result that $\dim(C_w)=|N(w)|=|N^-(w)|$ gives a combinatorial formula for the dimension of each Schubert cell.  Using \cite[Lemma 5.1]{P}, if $U\cdot x_J=\{ g\cdot x_J: g\in U \}$ then
\begin{eqnarray}\label{eq: U-action}
U\cdot x_J=x_J+\V; \;\; \mbox{  where  } \V=\bigoplus_{\gamma\in \Phi^+-J} \fg_{\gamma}.
\end{eqnarray}
Applying \cite[Proposition 3.7]{P} to equation \eqref{eq: U-action} we get the following dimension formula for the intersection $C_w\cap \B(x_J,H)$. 

%The follow Lemma is a special case of Lemma 5.1 from \cite{P}.  
% 
%\begin{lem}\label{lemma1}  The Adjoint $U$-orbit of $x_J$ in $\fg$ is,
%\[
%U\cdot x_J=x_J+\V; \;\; \mbox{  where  } \V=\bigoplus_{\gamma\in \Phi^+-J} \fg_{\gamma}
%\]
%and $U\cdot x_J:=\{Ad(u)(x_J) : u\in U \}$.
%\end{lem}

%Using Proposition 3.7 in \cite{P} we get the following dimension formula for the cell $C_w\cap \B(x_J,H)$. 

\begin{lem}\label{lemma2} For all $w\in W$, $C_w\cap \B(x_J,H)\neq \emptyset$ if and only if $w^{-1}(J)\subseteq \Phi_H$.  If $C_w\cap \B(x_J,H)\neq \emptyset$ then $\dim(C_w\cap \B(x_J,H)) = |N^-(w) \cap \Phi_H^-|$.  
\end{lem}
\begin{proof}  By \cite[Proposition 3.7]{P}, $C_w\cap \B(x_J,H)$ is nonempty if and only if $w^{-1}\cdot x_J \in H$.  Since
\[
w^{-1}\cdot x_J=w^{-1}\cdot n_J +w^{-1}\cdot s_J=\sum_{\alpha\in J}E_{w^{-1}(\alpha)} +w^{-1}\cdot s_J	
\]
and $w^{-1}\cdot s_J\in \fh\subseteq H$ for all $w$, we see that $w^{-1}\cdot x_J\in H$ if and only if $w^{-1}(\alpha)\in \Phi_H$ for all $\alpha\in J$ or equivalently $w^{-1}(J)\subset \Phi_H$.  Finally, applying the dimension formula from \cite[Proposition 3.7]{P} and simplifying we have,
\begin{eqnarray*}
\dim(C_w\cap \B(x_J,H)) &=&  |N(w^{-1})|- \dim \V/(\V \cap w\cdot H) \\
	&=& |N(w^{-1})|-|\{ \gamma\in \Phi^+-J : w^{-1}(\gamma)\notin \Phi_H \}|\\
	&=& |N(w^{-1})|-|\{ \gamma\in N(w^{-1}) : w^{-1}(\gamma) \notin \Phi_H^- \}|\\
	&=& |\{ \gamma\in N(w^{-1}): w^{-1}(\gamma)\in \Phi_H^- \}|\\
	&=&|N(w^{-1})\cap w(\Phi_H^-)|\\
	&=& |N^-(w)\cap \Phi_H^-|.
\end{eqnarray*}
The last equality follows from the fact that $w^{-1}(N(w^{-1}))=N^-(w)$ for all $w\in W$.
\end{proof}

To calculate the Betti numbers of $\B(x_J,H)$, we will use the fact that regular Hessenberg varieties are paved by affines.  The following is a special case of \cite[Theorem 5.4]{P}.

\begin{lem}\label{lemma3} Let $J\subset \Delta$ and $x_J$ be a regular element as in equation \eqref{regelement}.  Then the regular Hessenberg variety $\B(x_J,H)$ is paved by affines.  This paving is given by the intersections $\overline{C_w}\cap \B(x_J,H)$, and the $i$-th Betti number of $\B(x_J,H)$ is 
\[
\beta_i(J):=|\{w\in W: w^{-1}(J)\subseteq \Phi_H,\,\, |N^-(w)\cap \Phi_H^-|=i \}|.
\]
\end{lem}

%Throughout this paper, we'll assume that $-\Delta\subseteq \Phi_H$ which implies that the corresponding regular semisimple Hessenberg variety is connected (see \cite[Appendix 1]{AT}).  

%This results in no loss of generality for studying the representations of the Weyl group, since Teff has shown in \cite{T} that in this case the representation is given by the induced representation of the action on each component.   

%More generally, the proof methods in \cite{P2}, Theorem 3.4 give the following result.
%
%\begin{lem} Suppose there exists $-\alpha\in -\Delta$ such that $-\alpha\notin \Phi_H^-$.  Then there exists $w\in W$ such that 
%\end{lem}

\begin{cor}\label{dimension} The dimension of the regular Hessenberg variety $\B(x_J, H)$ is $m_H$.
\end{cor} 
\begin{proof}  Let $\Delta_H =\{ \alpha\in \Delta : -\alpha\in \Phi_H^- \}$ and $W_H$ denote the subgroup of $W$ generated by the simple reflections $s_{\alpha}$, $\alpha\in \Delta_H$.  If $\Phi_{\Delta_H}$ is the root system generated by $\Delta_H$ with negative roots $\Phi_{\Delta_H}^-$, then $\Phi_H^-\subseteq \Phi_{\Delta_H}^-$ since $H$ is a $\fb$-submodule of $\fg$.  Let $w_H\in W$ denote the longest element of $W_H$.  Notice that $N^-(w_H)= \Phi_{\Delta_H}^-$ and that $w_H^{-1}(\Delta) \cap \Phi^- = -\Delta_H \subseteq \Phi_H^-$ so $C_{w_H}\cap \B(x_J,H)\neq \emptyset$.  For all $w\in W$ such that $w^{-1}(J)\subseteq \Phi_H$ we have,
\begin{eqnarray*}
	\dim(C_w\cap \B(n_J,H)) &=& |N^-(w)\cap \Phi_H^-|\\ 
					&=& |N^-(w)\cap \Phi_H^-\cap \Phi_{\Delta_H}^-|\\
					&\leq& |\Phi_H^-\cap \Phi_{\Delta_H}^-|\\
					&=& \dim(C_{w_H} \cap \B(x_J,H)).
\end{eqnarray*}
Therefore $\overline{C_{w_H}}\cap \B(x_J,H)$ is a maximum dimensional cell in the paving given in Lemma \ref{lemma3} and so $\dim(\B(x_J,H))=\dim(C_{w_H}\cap \B(x_J,H))=|\Phi_H^-\cap \Phi_{\Delta_H}^-|= m_H$.
\end{proof}

%%%%%%%%%%%%%%%%%%%%%%%%%%%%%%%%%%%%%%%%%%%%

\section{Results}

To show that the Betti numbers of the regular Hessenberg variety $\B(x_J, H)$ are palindromic we study the sets $N^-(w)\cap \Phi_H^-$ in greater detail.  This study has already been initiated by Sommers and Tymoczko in \cite{ST}.  Following their terminology, we give the following definition.

\begin{defn}  Given a subset $S\subset \Phi_H^-$, we say $S$ is $\Phi_H^-$-closed if for all $\alpha, \beta\in S$ such that $\alpha+\beta\in \Phi_H^-$, then $\alpha+\beta\in S$ as well.  Given such a subset $S\subset \Phi_H^-$, we say that $S$ is of Weyl type if both $S$ and $S^c:=\Phi_H^- - S$ are $\Phi_H^-$-closed.  Denote the set of all subsets of Weyl type in $\Phi_H^-$ by $\W^H$.
\end{defn}
 
When $H=\fg$, these subsets are analogous to Weyl group elements since $S\subseteq \Phi^-$ is of Weyl type if and only if $S=N^-(w)$ for a unique $w\in W$ (see \cite[Proposition 5.10]{Ko}).   In addition, 
\begin{eqnarray}
S^c=\Phi^- - N^-(w)=N^-(w_0w),
\end{eqnarray}
where $w_0$ denotes the longest element of $W$.  The following result characterizes subsets of Weyl type in $\Phi_H^-$ and is a combination of \cite[Proposition 6.1]{ST} and \cite[Proposition 6.3]{ST}.

\begin{prop}\label{Weyl-type}  Let $S\in \W^H$ be a subset of Weyl type. 
\begin{enumerate}
\item There exists $w\in W$ such that $S=N^-(w)\cap \Phi_H^-$, and $S$ is a subset of Weyl type in $\Phi_H^-$ if and only if it is of this form.  
\item There exists a unique $w\in W$ satisfying both $S=N^-(w)\cap \Phi_H^-$ and $w^{-1}(\Delta)\subseteq \Phi_H$.
\end{enumerate}
\end{prop}

The combinatorial results of Proposition \ref{Weyl-type} give a proof of the fact that the Betti numbers of regular nilpotent Hessenberg varieties are palindromic, which is part of an unpublished result of Peterson (see \cite[Theorem 3]{BC}).  

%The proof showcases the effectiveness of using subsets of roots of Weyl-type.

\begin{prop}\label{regular nilpotent}  For all $w\in W$ such that $w^{-1}(\Delta)\subseteq \Phi_H$ there exists a unique $\bar{w}\in W$ such that $\bar{w}^{-1}(\Delta)\subseteq \Phi_H$ and $N^-(\bar{w})\cap \Phi_H^- = (N^-(w)\cap \Phi_H^-)^c$. In particular, the Betti numbers of regular nilpotent Hessenberg varieties are palindromic.
\end{prop}
\begin{proof}  Every regular nilpotent element of $\fg$ is conjugate to an element of the form $n=\sum_{\alpha\in \Delta} E_{\alpha}$ so without loss of generality, we assume that $n$ is of this form.  Consider $S=N^-(w)\cap \Phi_H^-$, which is a subset of Weyl type by Proposition \ref{Weyl-type}.  This implies that $S^c$ is also of Weyl type, and Proposition \ref{Weyl-type} also guarantees the existence of a unique $\bar{w}$ such that $\bar{w}^{-1}(\Delta)\subseteq \Phi_H$ and $N^-(\bar{w})\cap\Phi_H^- =S^c$ as we wanted.  By Lemma \ref{lemma3} the $i$-th Betti number of $\B(n,H)$ is
\[
\beta_i(\Delta)=|\{w\in W : w^{-1}(\Delta) \subseteq \Phi_H ,\,\, |N^-(w)\cap \Phi_H^-| = i \}|.
\]  
Since $|S^c|=|\Phi_H^-|-|S|=m_H - |S|$, we conclude that the Betti numbers of $\B(n,H)$ are palindromic.
\end{proof}

\begin{rem}  Note that regular nilpotent Hessenberg varieties can be singular.  The singular locus for $\B(n, H)$ when $\Phi_H^- = -\Delta$ and $n$ is a regular nilpotent element is described by Insko and Yong in \cite{IY}.
\end{rem}

For each $S\in \mathcal{W}^H$ and $J\subseteq \Delta$ set 
\[
	W(J,S)=\{w\in W: w^{-1}(J)\subseteq \Phi_H,\, S=N^-(w)\cap \Phi_H^-\}.
\]  
By Proposition \ref{Weyl-type}, this set is nonempty for all nonempty $S\in \W^H$ and contains a single element if $J=\Delta$.  By Lemma \ref{lemma3}, 
\begin{eqnarray}\label{Bettieq}
		\beta_i(J)=\sum_{S\in \W^H,\, |S|=i} |W(J,S)|.
\end{eqnarray} 
The following Lemma catalogues some well-known results on coset representatives.  
%These results can be found in \cite[\S 5]{Ko}.  

\begin{lem}\label{cosetreps}  Let $J\subseteq \Delta$ and $W_J$ be the corresponding subgroup of the Weyl group generated by the simple reflections $s_{\alpha}$ for $\alpha\in J$.  The set of shortest coset representatives of $W_J\backslash W$ is
\[
	W^J=\{ v\in W : N(v^{-1}) \subseteq \Phi^+-\Phi_J^+ \}.
\]
Each $w\in W$ can be written uniquely as $w=yv$ with $y\in W_J$ and $v\in W^J$, and this decomposition has the property that $\ell(w)=\ell(y)+\ell(v)$ and therefore
\begin{eqnarray}\label{root decomp}
	N^-(yv)=N^-(v)\sqcup v^{-1}N^-(y).
\end{eqnarray}  
\end{lem}

It's also a well-known fact that $W_J$ normalizes $\Phi^+-\Phi_J^+$ (and therefore $\Phi^--\Phi_J^-$) so for any $y\in W_J$,
\begin{eqnarray}\label{normalize}
	y(\Phi^- - \Phi_J^-)=\Phi^- - \Phi_J^-.
\end{eqnarray}
Recall that for any $J\subseteq \Delta$, $\fm_J$ denotes the Levi subalgebra of $\fg$ corresponding to the subroot system $\Phi_J$ in $\Phi$.  The following result is a combination of \cite[Proposition 5.2]{P} and \cite[Corollary 5.8(2)]{P}.    
 
\begin{lem}\label{lemma4}  For each subset $J\subseteq \Delta$ and $v\in W^J$, 
\[
	H_v:=v\cdot H\cap \fm_J 
\] 
is a Hessenberg space in $\fm_J$ with respect to $\fb_J:= \fb\cap \fm_J$. Further, let $w=yv$ with $y\in W_J$ and $v\in W^J$.  Then $C_w\cap \B(x_J,H)\neq \emptyset$ (i.e., $w^{-1}(J)\subseteq \Phi_H$) if and only if $C_y\cap \B(n_J, H_v) \neq \emptyset$ ($i.e., y^{-1}(J)\subseteq \Phi_{H_v}$).
\end{lem}  

\begin{rem}\label{rootspace-decomp} Lemma \ref{lemma4} allows us to relate subsets of Weyl type in $\Phi_H^-$ to subsets of Weyl type in $\Phi_{H_v}^-$.  Applying Lemma \ref{cosetreps} and the fact that $N^-(y)\subseteq \Phi_J$, we get that
\begin{eqnarray}\label{Hess decomp}
	N^-(yv)\cap\Phi_H^- &=& (N^-(v)\cap \Phi_H^-) \sqcup (v^{-1}N^-(y)\cap \Phi_H^-)\nonumber\\ 
				&=&  (N^-(v)\cap \Phi_H^-) \sqcup v^{-1}(N^-(y)\cap v\cdot \Phi_H^-)\\
				&=&(N^-(v)\cap \Phi_H^-) \sqcup v^{-1}(N^-(y)\cap \Phi_{H_v}^- )\nonumber.
\end{eqnarray}   
\end{rem}

\begin{prop}\label{cosetsneq}  Let $w$ and  $x$ be distinct elements of $W$ with $w=yv$ and $x=zu$ where $y,z\in W_J$ and $u,v\in W^J$.  If $w, x \in W(J,S)$ for $S\in \mathcal{W}^H$ then $v\neq u$.
\end{prop}
\begin{proof}  Seeking a contradiction, suppose $v=u$ and let's denote this element by $v$.  Note that as we assume $w$ and $x$ are distinct it must be the case that $y\neq z$.  Since $w,x\in W(J,S)$ we have that 
\[
	N^-(w)\cap \Phi_H^- =S= N^-(x)\cap \Phi_H^-
\]
and $w^{-1}(J), x^{-1}(J)\subseteq \Phi_H$.  Using equation \eqref{Hess decomp} and Lemma \ref{lemma4} we see that
\[
	(N^-(v)\cap \Phi_H^-) \sqcup v^{-1}(N^-(y)\cap \Phi_{H_v}^-) = (N^-(v)\cap \Phi_H^-) \sqcup v^{-1}(N^-(z)\cap \Phi_{H_v}^-)
\]
and $y^{-1}(J), z^{-1}(J)\subseteq \Phi_{H_v}$.  In particular it must be the case that,
\[
N^-(y)\cap \Phi^-_{H_v} = N^-(z)\cap \Phi^-_{H_v}.
\]
But now $y$ and $z$ are distinct elements of $W_J$ such that $y^{-1}(J), z^{-1}(J) \subseteq \Phi_{H_v}$ and $N^-(y)\cap \Phi^-_{H_v} = N^-(z)\cap \Phi^-_{H_v}$, contradicting the results of Proposition \ref{Weyl-type}(2) which state that such an element must be unique.
\end{proof}

This proposition shows that the elements of $W(J,S)$ must correspond to different cosets of $W_J\backslash W$.  To show that the Betti numbers of $\B(x_J, H)$ are palindromic, it is enough to show that the elements of $W(J,S)$ are in bijective correspondence with the elements of $W(J, S^c)$ by equation \eqref{Bettieq}.  Given $w\in W(J,S)$, a natural candidate for the corresponding element of $W(J,S^c)$ is $w_0w$ since $N^-(w_0w)\cap\Phi_H^- = S^c$.  

%, which in turn would imply that
%	\begin{eqnarray*} 
%		\beta_i(J)&=& \sum_{S\in \W^H,\, |S|=i} |W(J,S)|\\
%			&=& \sum_{S \in \W^H, |S|=i} |W(J, S^c)|\\
%			&=& \sum_{R\in W^H, |R|=m_H-i} |W(J,R)|\\
%			&=& \beta_{m_H-i}(J)
%	\end{eqnarray*}
%since every $R\in \W^H$ such that $|R|=m_H-i$ must be of the form $R=S^c$ for some $S\in \W^H$ such that $|S|=i$.  Suppose $w\in W(J,S)$. A natural candidate for a corresponding $\bar{w}\in W(J,S^c)$ is $w_0w$ since then $N^-(w_0w)=\Phi^- - N^-(w)$ and therefore $N^-(w_0w)\cap\Phi_H^- = S^c$ as a subset of $\Phi_H^-$.  

The difficulty is as follows.  Although $w^{-1}(J)\subseteq \Phi_H^-$ it need not be the case that $(w_0w)^{-1}(J)\subseteq \Phi_H^-$, so $w_0w$ may not be an element of $W(J,S^c)$.  Furthermore, if $w$ and $x$ are distinct elements of $W(J,S)$, then they belong to different cosets of $W_J\backslash W$.  This need not be the case for $w_0w$ and $w_0x$.  To address this difficulty, we instead consider the subset $K=-w_0(J)\subseteq \Delta$ of simple roots.

\begin{prop}\label{J and K cosets}  Let $w\in W$ and write $w=y_Jv_J$ with $y_J\in W_J$ and $v_J\in W^J$.  Similarly, write $w_0w=y_K v_K$ for $y_K \in W_K$ and $v_K\in W^K$.  Then, 
\[
	N(v_K)=v_J^{-1}(\Phi^+-\Phi_J^+)\cap \Phi^+  \mbox{  and  } N(v_J)=v_K^{-1}(\Phi^+- \Phi_K^+) \cap \Phi^+
\]
and the coset of $w_0w$ in $W_K\backslash W$ is uniquely determined by the coset of $w$ in $W_J \backslash W$.
\end{prop}
\begin{proof}  We'll show that $N(v_J)=v_K^{-1}(\Phi^+- \Phi_K^+)\cap \Phi^+$.  The fact that $N(v_K)=v_J^{-1}(\Phi^+-\Phi_J^+)\cap \Phi^+$ will follow from the same argument with the roles of $K$ and $J$ reversed.  By definition, $\gamma\in N(v_J)$ if and only if $\gamma\in \Phi^+$ such that $v_J(\gamma)\in \Phi^-$.  This is the case if and only if $v_J(\gamma)\in N^-(v_J^{-1})$, and it follows that $v_J(\gamma)\in \Phi^- - \Phi_J^-$ by Lemma \ref{cosetreps}.  We apply equation \eqref{normalize} and the fact that $w_0(\Phi^-)=\Phi^+$ and $w_0(\Phi_J^-)=\Phi_K^+$ to obtain,
\begin{eqnarray*}
	v_J(\gamma)\in \Phi^- -\Phi_J^- &\Leftrightarrow &  y_Jv_J(\gamma)\in y_J(\Phi^--\Phi_J^-)=\Phi^--\Phi_J^-\\
				&\Leftrightarrow& y_Kv_K(\gamma)=w_0y_Jv_J(\gamma)\in w_0(\Phi^- -\Phi_J^-)=\Phi^+-\Phi_K^+\\
				&\Leftrightarrow & v_K(\gamma)\in y_K^{-1}(\Phi^+ -\Phi_K^+)=\Phi^+ - \Phi_K^+\\
				&\Leftrightarrow& \gamma\in v_K^{-1}(\Phi^+-\Phi_K^+).
\end{eqnarray*}
Therefore $\gamma\in N(v_J)$ if and only if $\gamma\in \Phi^+$ and $\gamma\in v_K^{-1}(\Phi^+-\Phi_K^+)$, or equivalently $\gamma\in v_K^{-1}(\Phi^+-\Phi_K^+)\cap \Phi^+$.

To show that this property uniquely determines the cosets, recall that any $v\in W$ is uniquely determined by the set $N(v)$. Suppose $w=y_Jv_J$ and $x=z_Ju_J$ for $y_J, z_J\in W_J$ and $v_J, u_J\in W^J$, and that $v_J\neq u_J$.  If $w_0w=y_Kv_K$ and $w_0x= z_Ku_K$ for $y_K, z_K\in W_K$ and $v_K, u_K\in W^K$, we claim that $v_K\neq u_K$.  Indeed, if not then $v_K=u_K$ and
	\[
		N(v_J)=v_K^{-1}(\Phi^+-\Phi_K^+)\cap \Phi^+ =  u_K^{-1}(\Phi^+-\Phi_K^+)\cap \Phi^+=N(u_J)
	\]
and we conclude that $v_J=u_J$, a contradiction.
\end{proof}

We can finally prove our main theorem.

\begin{thm}\label{main result}  The Betti numbers of a regular Hessenberg variety are palindromic.
\end{thm}
\begin{proof}  By equation \eqref{regelement}, it suffices to prove this statement for $\B(x_J,H)$.  We'll prove that $\beta_i(J)=\beta_{m_H-i}(K)$ where $K=-w_0(J)$.  The claim will then follow from the fact that $\B(x_J,H)$ and $\B(x_K, H)$ are isomorphic (since $x_J$ and $x_K$ are conjugate) so their Betti numbers are equal.  

Fix $S\in \W^H$ such that $|S|=i$.  To prove that $\beta_i(J)=\beta_{m_H-i}(K)$, by equation \eqref{Bettieq} it suffices to give a bijection between $W(J,S)$ and $W(K,S^c)$ for any such $S\in \W^H$.  Let $w\in W(J,S)$ and write $w_0w=yv$ with $y\in W_K$ and $v\in W^K$.  As noted above, $w_0w$ is not a priori an element of $W(K,S^c)$ although it does satisfy 
\[
	N^-(w_0w)\cap \Phi_H^- = (\Phi^- - N^-(w)) \cap \Phi_H^- = S^c.  
\]
By Lemma \ref{lemma4} it must be the case that $H_v=v\cdot H\cap \fm_K$ is a Hessenberg space of $\fm_K$.  Proposition \ref{Weyl-type} implies that $N^-(y)\cap \Phi_{H_v}$ is a subset of Weyl-type in $\Phi_{H_v}$ and that there exists a unique $\bar{y}\in W_K$ such that $\bar{y}^{-1}(K)\subseteq \Phi_{H_v}$ and $N^-(\bar{y})\cap \Phi_{H_v}=N^-(y)\cap \Phi_{H_v}$. Set $\bar{w}=\bar{y}v$.  Then $\bar{w}^{-1}(K)\subseteq \Phi_H$ by Lemma \ref{lemma4} and
\begin{eqnarray*}
N^-(\bar{w})\cap \Phi_H^-&=& (N^-(v)\cap \Phi_H^-) \sqcup v^{-1}(N^-(\bar{y})\cap \Phi_{H_v}^-)\\
				&=& (N^-(v)\cap \Phi_H^-) \sqcup v^{-1}(N^-(y)\cap \Phi_{H_v}^-)\\
				&=& N^-(w_0w)\cap \Phi_H^-\\
				&=& S^c
\end{eqnarray*}
so $\bar{w}\in W(K,S^c)$.  In this way, we obtain a map $W(J,S)\to W(K, S^c)$ given by $w\mapsto \bar{w}$.

We now need to show that this association is a bijection.  By Proposition \ref{cosetsneq} two distinct elements $w,x\in W(J,S)$ have distinct cosets in $W_J\backslash W$, so if $w= y_J v_J$ and $x=z_J u_J$ where $y_J,z_J\in W_J$ and $v_J, u_J\in W^J$, then $v_J\neq u_J$.  If $w_0w=y_Kv_K$ and $w_0x=z_Ku_K$ where $y_K, z_K\in W_K$ and $v_K, u_K\in W^K$ then $v_K\neq u_K$ by Proposition \ref{J and K cosets}.  Since $v_K\neq u_K$, $\bar{w}$ and $\bar{x}$ are distinct elements of $W(K,S^c)$ and we see that the map is injective. Repeating the argument above with the roles of $K$ and $J$ exchanged gives an injective map $W(K,S^c) \to W(J,S)$ so $|W(K,S^c)| \leq |W(J,S)|$ and our original map must therefore be a bijection.  

%Finally, since $|W(J,S)|=|W(K,S^c)|$ we get,
%\begin{eqnarray*}
%	\beta_i(J) &=& \sum_{S\in \W^H,\, |S|=i} |W(J,S)|\\
%			&=& \sum_{S\in \W^H,\, |S|=i} |W(K,S^c)|\\ 
%			&=& \sum_{R\in \W^H, |R|=m_H-i} |W(K,R)|\\
%			&=& \beta_{m_H-i}(K)
%\end{eqnarray*}
%where the last equality follows from the fact that every $R\in \W^H$ such that $|R|=m_H-i$ must be of the form $R=S^c$ for some $S\in \W^H$ such that $|S|=i$ by definition of Weyl type subsets.
\end{proof}

We emphasize that the result of Theorem \ref{main result} is notable since regular Hessenberg varieties are not in general smooth, and leads us to speculate that these varieties are rationally smooth.  

We further remark that combining Corollary \ref{dimension} with the methods used by Anderson and Tymoczko in \cite[Lemma 7.1]{AT} yields the following statement.  This seems well know by experts in the field, but we state it here for future reference.  

%The condition that $-\Delta\subseteq \Phi_H$ implies that the corresponding regular semisimple Hessenberg variety is connected (see \cite[Appendix 1]{AT}), and therefore the regular Hessenberg variety is as well.

\begin{cor}  If $-\Delta \subseteq \Phi_H^-$ then the regular Hessenberg variety $\B(x_J,H)$ is irreducible.  
\end{cor}
\begin{proof}  The condition that $-\Delta\subseteq \Phi_H^-$ implies that the element $w_H$ from the proof of Corollary \ref{dimension} is $w_0$.  In this case, $\overline{C_{w_0}}\cap \B(x_J, H)$ is the unique maximal dimensional cell in the paving given in Lemma \ref{lemma3} and therefore $\B(x_J,H)$ is connected by Theorem \ref{main result}.  Finally, the proof of \cite[Lemma 7.1]{AT} implies that $\B(x_J,H)$ is pure-dimensional of dimension $m_H$.  It immediately follows that $\B(x_J,H)= \overline{C_{w_0}}\cap \B(x_J,H)$.
\end{proof}

We hope that the results of this paper aid understanding of regular Hessenberg varieties and lead to further results describing their geometry and the Weyl group representation discussed in the introduction.

\end{document}